\documentclass[10pt]{IEEEtran}

\usepackage{xspace,amsmath,amssymb,epsfig,subfigure,syntonly,amsthm}
\usepackage{epsfig,amsmath,color}
\usepackage{xspace,syntonly}
\usepackage{algorithm,algorithmic}

\newcommand{\utwi}[1]{\mbox{\boldmath $#1$}}

\newcommand{\Prob}{{\mathrm{Pr}}}

\def \pnmax {{p_n^\mathrm{max}}}
\def \Pnmax {{P_n^\mathrm{max}}}
\def \iRmax {{\iota_R^\mathrm{max}}}
\def \inRmax {{\iota_{R,n}^\mathrm{max}}}

\def \ntoi {n \rightarrow i}
\def \ntoj {n \rightarrow j}
\def \jton {j \rightarrow n}
\def \iton {i \rightarrow n}

\newcommand{\cL}{{\cal{L}}}
\newcommand{\cN}{{\cal N}}
\newcommand{\cP}{{\cal P}}

\newcommand{\cA}{{\cal A}}

\newcommand{\cC}{{\cal C}}

\newcommand{\cI}{{\cal I}}
\newcommand{\cB}{{\cal B}}
\newcommand{\cU}{{\cal U}}

\newcommand{\bq}{{\bf q}}

\newcommand{\bx}{{\bf x}}

\newcommand{\bD}{{\bf D}}

\newcommand{\bvartheta}{{\utwi{\vartheta}}}

\newtheorem{proposition}{Proposition}

\begin{document}

\IEEEoverridecommandlockouts

\title{Statistical Routing for Multihop \\ Wireless Cognitive Networks}

\author{Emiliano Dall'Anese, \emph{Member, IEEE}, and Georgios B. Giannakis, \emph{Fellow, IEEE}

\thanks{Manuscript submitted on December 23, 2011; revised May 14, 2012; accepted July 2, 2012. This work was supported by the QNRF grant NPRP 09-341-2-128. Part of the paper appeared in the \emph{Proc. of Intl. Conf. on Acoust., Speech, and Signal Proc.}, Kyoto, Japan, March 2012.
}
\thanks{\protect\rule{0pt}{0.5cm}%
The authors are with the Department of Electrical and Computer Engineering, University of Minnesota, 200 Union Street SE, %
Minneapolis, MN 55455, USA. Tel/fax: +1(612)624-9510/625-2002. %
E-mails: {\tt \{emiliano, georgios\}@umn.edu}
}
}

\markboth{IEEE JOURNAL ON SELECTED AREAS IN COMMUNICATIONS (TO APPEAR)}%
{First Author \MakeLowercase{\textit{et al.}}: Title}

\maketitle

\begin{abstract}
To account for the randomness of propagation channels and interference levels in hierarchical spectrum sharing, a novel approach to multihop routing is introduced for cognitive random access networks, whereby packets are randomly routed according to outage probabilities.  Leveraging channel and interference level statistics, the resultant cross-layer optimization framework provides optimal routes, transmission probabilities, and transmit-powers, thus enabling cognizant adaptation of routing, medium access, and physical layer parameters to the propagation environment. The associated optimization problem is non-convex, and hence hard to solve in general. Nevertheless, a successive convex approximation approach is adopted to efficiently find a Karush-Kuhn-Tucker solution. Augmented Lagrangian and primal decomposition methods are employed to develop a distributed algorithm, which also lends itself to online implementation. Enticingly, the fresh look advocated here permeates benefits also to conventional multihop wireless networks in the presence of channel uncertainty. 
\end{abstract}


\IEEEpeerreviewmaketitle

\section{Introduction}
\label{sec:Introduction}

Resembling traditional routing protocols for wired networks, their counterparts for wireless networking generally utilize optimization tools such as shortest path routing to find optimal route(s) based on the network connectivity graph abstraction~\cite{Silvester93}. Early on, links among nodes were quantified based on a disk model capturing only distance-based deterministic losses. Upon recognizing the inadequacy of disk models for the broadcast wireless interface~\cite{Haenggi05}, a weighted graph accommodating more sophisticated performance metrics was adopted; see e.g.,~\cite{DeCouto03},~\cite{ZLZhang07}, and the stochastic routing approach in~\cite{Ribeiro-infocom07}, where link weights capturing packet delivery probabilities were exploited to develop optimal routing schemes. These schemes are particularly attractive for energy-limited nodes, primarily because the resulting routing strategies promote links with higher reliability, thus decreasing the number of packet lost due to fading~\cite{Ephremides-twc02}.  

In a hierarchical access setting, interference levels can not be acquired accurately due to the lack of explicit inter-system cooperation~\cite{ZhSa07}. As a result, random shadowing and small-scale fading effects, along with dynamically changing activities of licensed users, accentuate the uncertain nature of wireless cognitive radio (CR) links. The effects of random interference on CR links from primary user (PU) transmitters is called upon in~\cite{Khalife08}, where source-to-destination paths that are most likely to meet prescribed end-to-end requirements are found based on predicted link capacities. 
Leveraging the situational-awareness provided by spectrum occupancy detection schemes, a graph whose link weights reflect the amount of spectral resources available per CR-to-CR link is employed in~\cite{Xin08}, where optimal routes are obtained via Dijsktra or Bellman Ford-like algorithms. A two-phase approach is proposed in~\cite{Pefkianakis08}, where nodes in the network first obtain an expected route cost and a set of candidate forwarding nodes, and then route traffic across paths with higher
spectrum availability. In~\cite{abbagnale-Secon}, the average link availability is invoked to develop a routing scheme that avoids network zones with unstable CR connectivity. Link availability in~\cite{abbagnale-Secon} is computed in a probabilistic sense based only on the statistics of primary user (PU) activities. A PU coverage map supplied by sensing schemes is employed in~\cite{Chowdhury11} to identify spectrum opportunities in space, and devise routing strategies supporting multiple classes of CR quality-of-service (QoS) demands. 

The aforementioned works offer valuable insights on route formation and management based on the average availability of CR links, and predicted link capacities. However, in a hierarchical access setup, link capacities are unknown and may change abruptly because of time-varying PU activity patterns, dynamic shadowing, and diverse QoS constraints. In this context, a cross-layer design approach to obtain \emph{both} optimal routes \emph{and} physical and medium access parameters that dictate the packet forwarding capabilities is therefore well motivated. To this end, the present paper exploits propagation channel statistics to develop a statistical routing approach whereby nodes not only compute optimal routes, but also optimal link reliabilities by controlling transmit-powers and medium access control (MAC) parameters. The novel approach accounts explicitly for the randomness of propagation and the medium access interface, to allow spectrum-cognizant routing of data packets, while enforcing PU interference protection (Section~\ref{sec:statistical}).  

In spite of the non-convexity of the associated cross-layer optimization problem, a successive convex approximation is pursued to find a Karush-Kuhn-Tucker (KKT) solution efficiently (Section~\ref{sec:tractable}). Enticingly, feasibility guarantees offered by the successive convex approximation algorithm naturally suggest an online implementation of the algorithm whereby nodes do not necessarily wait for the successive convex approximation iterations to converge, but rather use network parameters as they become available. 

However, the communication overhead incurred to acquire channel statistics at a central node, and subsequently disseminate optimal network parameters can become prohibitive as the network size increases. To alleviate such a message-passing burden, and address scalability and robustness concerns, a distributed algorithm is also developed by invoking the alternating direction method of multipliers and the primal decomposition method (Section~\ref{sec:distributed}). Finally, suitable conditions are established to ensure that packets are eventually delivered to their destination when routes, medium access and physical layer parameters are regularly updated to track channel statistics and topology dynamics (Section~\ref{sec:deliverability}).

\subsection{Preliminaries and problem formulation}
\label{sec:preliminaries}

Consider a wireless CR network with $N$ nodes $\{U_n\}_{n=1}^{N}$ 
sharing spectral resources with an incumbent PU system~\cite{ZhSa07}. 
Leveraging the spectrum awareness provided by spatio-temporal sensing
schemes~\cite{KDGjstsp11,DBGphycom}, CRs collaborate in routing data packets to a sink node $U_{N+1}$, while respecting the PU-CR hierarchy. The CR network is modeled as a digraph to account for the
possible lack of link bi-directionality.
The dynamic and stochastic nature of the CR propagation ambience, along with the possibly minimal amount of topological information motivate consideration of random medium access, as well as stochastic routing strategies~\cite{ZLZhang07, Ribeiro-infocom07,Lott06}. 
In this context, a CR node $U_n$ transmits with probability $\mu_n \in [0,1]$, and decides whether to route packets toward a neighboring node $U_i$ 
with probability $t_{n \rightarrow i} \in [0,1]$ per time slot. As packets are forwarded to neighboring nodes according to probability mass functions, it holds that 
$\sum_{i \neq n} t_{\ntoi} = 1$, for all $n = 1,\ldots,N$.

Communication of data packets over a wireless network depends not only on transmission and forwarding decisions, but also on the intended link reliability. In case of unsuccessful packet decoding due to fading- or interference-induced link outages~\cite{Haenggi05}, a packet not eventually routed by $U_n$ will remain in $U_n$'s queue, and its transmission will be re-attempted in a subsequent time slot (possibly to a different neighboring CR). To capture channel- and interference-induced sources of uncertainty, let $r_{n \rightarrow i} \in (0,1]$ denote the probability that a packet transmitted from node $U_n$ is correctly decoded (and thus successfully received) by $U_i$.

Assuming that link reliabilities $\{r_{n \rightarrow i}\}$ are \emph{known} by, e.g., computing the packed error rate of preceding sessions, a stochastic routing framework for maximizing users' exogenous rates was introduced in~\cite{Ribeiro-infocom07}. However, because of the volatile CR channel characteristics, time-varying PU activity patterns, and diverse QoS constraints, $\{r_{n \rightarrow i}\}$ may change abruptly during the network operation. Hence, $\{r_{n \rightarrow i}\}$ may \emph{not} be known in advance. Building on first- and second-order statistics of the PU interference, as well as those of node-to-node channels, a statistical routing approach yielding optimal (i) routes, (ii) transmission probabilities, and (iii) transmit-powers is put forward in the ensuing section.

\section{Statistical Routing Framework}
\label{sec:statistical}

Data percolation through a wireless network is captured by the product packet delivery probabilities $\{t_{n \rightarrow i}\, r_{n \rightarrow i}\}$.
When random access is employed as MAC, it is common to consider a packet lost when collisions among CR transmissions occur. With $\cI_{ni}$ denoting the set of nodes whose transmissions interfere with link $U_n \rightarrow U_i$, the probability of collision-free packet transmission from $U_n$ to $U_i$ is given by $\prod_{j \in \cI_{ni}} (1-\mu_j)$. A widely-accepted criterion for successful packet reception is to require the signal-to-interference-plus-noise ratio (SINR) to stay above a certain threshold~\cite{Haenggi05,Ephremides-twc02}, which is generally determined by the receiver structure, modulation, and coding scheme. Let $g_{\ntoi}$ denote the channel gain between $U_n$ and $U_i$, modeling the effects of path loss, log-normal shadowing, and Nakagami-$m$ small-scale fading~\cite{Stu01}. Then, the SINR of link $U_n \rightarrow U_i$ can be expressed as
\begin{equation}
\gamma_{\ntoi} := \frac{p_n g_{\ntoi}}{\sigma_{i}^2 + \sum_{S = 1}^{N_S} \pi_{S,i} } 
\label{eq:sinr}
\end{equation} 
where $\sigma_i^2$ stands for the receiver noise power at $U_i$; $p_n \in (0,\pnmax]$ denotes the transmission power of $U_n$; and $\pi_{S,i}$ the interference perceived from PU transmitter $S = 1,\ldots, N_S$. Randomness of $\{\gamma_{\ntoi}\}$ in~\eqref{eq:sinr} emerges due to the shadowing and small-scale effects on the PU interference $\{\pi_{S,i}\}$. Furthermore, CR-to-CR gains $\{g_{\ntoi}\}$ may be known imperfectly because of insufficient time for channel training. Nonetheless, CR-to-CR and PU-to-CR deterministic path losses, and statistics of shadowing and small-scale fading can be acquired and used. To this end, it is useful to recall that the distribution of channel gains $\{g_{\ntoi}\}$ can be approximated as log-normal~\cite[Ch. 2]{Stu01},~\cite{DKGGPtwc11}. Furthermore, the Fenton-Wilkinson result~\cite{Fen60} asserts that the distribution of SINRs $\{\gamma_{\ntoi}\}$ in~\eqref{eq:sinr} can be well-approximated as log-normal too, with mean and variance expressed in terms of the first- and second-order moments of $\{g_{\ntoi}\}$ and $\{\pi_{S,i}\}$; 
see~\cite{DKGGPtwc11} for a detailed derivation. Consequently, $\Gamma_{\ntoi} := 10 \log_{10} \gamma_{\ntoi}$ will be approximately Gaussian distributed with mean $P_n + m_{\ntoi}$, where $P_n := 10 \log_{10} p_n$, and variance denoted by  
$\sigma_{\ntoi}^2$. The probability $r_{\ntoi}$ that a packet transmitted from $U_n$ is correctly received by $U_i$ can thus be expressed as
\begin{align}
r_{\ntoi} & = \prod_{j \in \cI_{ni}} (1-\mu_j) \Prob\{ \gamma_{\ntoi} > \bar  \gamma_{\ntoi} \} \nonumber \\
&\approx \prod_{j \in \cI_{ni}} (1-\mu_j) \, Q\left(\frac{\bar \Gamma_{\ntoi} - P_n - m_{\ntoi} }{\sigma_{\ntoi}}\right) 
\label{eq:rni}
\end{align}
where $Q(x) := \int_x^\infty \frac{1}{\sqrt{2 \pi}} e^{-\frac{x^2}{2}} \mathrm{d}x$ is the standard Gaussian tail function, $\bar \gamma_{\ntoi}$ is a prescribed SINR threshold, and $\bar \Gamma_{\ntoi} := 10 \log_{10} \bar \gamma_{\ntoi}$. 
Similar to~\cite{Shin04},  the main interest here is in the tail of the complementary cumulative density function (ccdf) of the SINR; in this case, the Fenton-Wilkinson method is known to provide accurate approximations for all the propagation scenarios of practical interest~\cite[Ch. 3]{Stu01},~\cite{DKGGPtwc11,Beaulieu94}.

\begin{figure}[t]
\begin{center}
\includegraphics[width=0.4\textwidth]{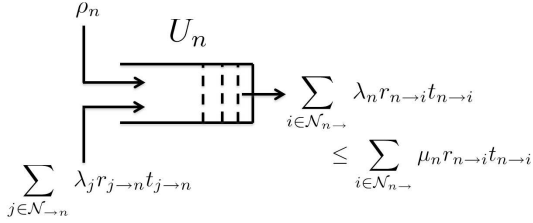}
\caption{Input and output flows at node $U_n$ under queue stability.}
\label{fig:queue}
\end{center}
\end{figure}

Using~\eqref{eq:rni}, the link reliabilities $\{r_{\ntoi}\}$ can be expressed in terms of the MAC variables $\{\mu_n\}$ and the physical layer quantities $\{P_n, m_{\ntoi},\sigma_{\ntoi}\}$. Therefore, with $\{m_{\ntoi},\sigma_{\ntoi}\}$ known parameters, the optimal routing strategy will be obtained by optimizing over $\{\mu_n, t_{\ntoi}, P_m\}$. The next step is to model exogenous data packet arrivals at $U_n$ from its application layer by a stationary stochastic process with average rate $\rho_n \in (0,1]$ per time slot. Suppose also that each CR node maintains a backlog to cache exogenous and endogenous\footnote{``Exogenous'' packets of a CR node are those generated from its application layer. On the other hand, ``endogenous''  packets refer to those received from the neighboring nodes of a CR node, and are to be routed by the network layer; see also Fig.~\ref{fig:queue}.} packets that have to be routed toward the destination $U_{N+1}$. Aggregate queue service rates depend on
the joint queue occupancy distribution. This results in a generally asymmetric system of interacting queues, whose stability region is challenging to analyze even for simple systems. Nevertheless, assuming as usual fully backlogged queues per node~\cite{Raotit88} yields a sufficient condition for queue stability that can be conveniently used as a constraint in rate-oriented routing optimization. In the advocated dominant system, users with empty queues transmit ``dummy'' packets, and consequently queue sizes are never smaller than those in the original system, if both systems start from the same initial condition. 

Let $\lambda_n$ denote the average aggregate rate of endogenous packet arrivals at $U_n$, which coincides with the rate of packet departures if queues are stable. Then, queue stability implies that $\{\rho_n\}_{n=1}^N$ and $\{\lambda_n\}_{n=1}^N$ 
abide by the flow conservation constraints~\cite{Ribeiro-infocom07,Georgiadis06} (cf. Fig.~\ref{fig:queue})
\begin{equation}
\rho_n = \lambda_n \sum_{j \in \cN_{n\rightarrow} } t_{n \rightarrow j} r_{n \rightarrow j} -  \sum_{i \in \cN_{\rightarrow n} } \lambda_i t_{i \rightarrow n} r_{i \rightarrow n}, \, \forall \, n 
\label{eq:consevation}
\end{equation}
where $\cN_{n \rightarrow} := \{j|r_{n \rightarrow j}>0 , j = 1,\ldots,N+1, j \neq n\}$ is the set of nodes that decode $U_n$'s transmissions with non-zero probability, and $\cN_{\rightarrow n} := \{i|r_{i \rightarrow n}>0 , i = 1,\ldots,N, i \neq n\}$ the set of nodes that route packets through $U_n$. 
For queue stability, Loynes' Theorem~\cite{Loynes} asserts that for stationary arrival and departure processes (the latter are stationary in the dominant system) a sufficient condition for stability is $\lambda_n < \mu_n$, for each CR $U_n$; and a necessary condition for stability is $\lambda_n \leq \mu_n$ (cf. Fig.~\ref{fig:queue}). 

To complete the formulation, consider $N_R$ actual or potential PU receivers, whose locations have been estimated via sensing~\cite{DKGtvt11}, and let $\iRmax$ denote the maximum average interference that can be tolerated by PU receiver $R$~\cite{ZhSa07,Zhang09}. Further, let $\cN_R \subseteq \{U_n\}_{n=1}^N$ be the (sub-)set of CR nodes located in the proximity of PU $R$ (not necessarily the entire CR network, as some CRs may be sufficiently far apart and do not interfere with PU $R$). Transmissions by CR $U_n$ undergo random shadowing and small-scale fading effects before arriving at close-by PU nodes. Approximate the channel gain $g_{n \rightarrow R}$ between CR $U_n$ and the PU $R$ as log-normal distributed~\cite{DKGGPtwc11}, and define a binary random variable $a_n \in \{0,1\}$, independent of $g_{n \rightarrow R}$, taking the value $1$ with probability $\mu_n$, and $0$ with probability $1-\mu_n$. Then, the average interference experienced at PU $R$ is given by ($\kappa := 0.1 \ln(10)$) 
\begin{align}
\iota_R & := \mathbb{E}\left\{\sum_{n \in \cN_R}a_n p_n g_{n \rightarrow R} \right\} \nonumber \\
& = \sum_{n=1}^N \mu_n e^{\kappa P_n + \kappa m_{n \rightarrow R} + \frac{\kappa^2}{2} \sigma^2_{n \rightarrow R}} \leq \iRmax .
\label{eq:interf}
\end{align}

Variables $\{P_n, \mu_n, \rho_n, \lambda_n\}$, and $\{t_{\ntoi}\}$ satisfying the constraints~\eqref{eq:consevation} and \eqref{eq:interf}
can be supported by the wireless CR network. It is certainly desirable
to design the network by selecting a feasible set of variables that
are optimal in some sense. To this end, consider a concave utility
$\cU_n(\rho_n)$, and a convex cost $\cC_n(P_n)$,
representing the reward of exogenous rate $\rho_n$ and the cost of power $P_n$ for node $U_n$, respectively~\cite{Marques_book}. Notice that $\rho_n$ is the average rate of packets  generated at the application layer of node $U_n$ to be eventually delivered to the sink $U_{N+1}$~\cite{Georgiadis06, Marques_book}; thus, $\rho_n$ represents an end-to-end performance metric. Capitalizing on the statistical description of SINRs and CR-to-PU channels, the \emph{statistical} routing problem is formulated as:
\begin{subequations}
%
%
\begin{align}
\textrm{(P1)} \max_{ \substack{\{P_n\}, \{\rho_n \geq 0\}, \\ \{\mu_n \geq 0\},  \{\lambda_n \geq 0\} \\ \{t_{\ntoi} \geq 0\}} } & \,\,\, \sum_{n=1}^N \cU_n (\rho_n) - \sum_{n=1}^N \cC_n(P_n) & \label{eq:sr_cost} \\
\textrm{subject to } & \nonumber \\
& \hspace{-2.5cm} \rho_n +  \sum_{i \in \cN_{\rightarrow n} } \lambda_i t_{i \rightarrow n} r_{i \rightarrow n} \leq \lambda_n \sum_{j \in \cN_{n \rightarrow} } t_{n \rightarrow j} r_{n \rightarrow j} \nonumber \\
& \hspace{2.3cm} \forall \,  n = 1,\ldots, N \label{eq:sr_flow} \\
\hspace{-2cm} \sum_{i \in \cN_{n\rightarrow}} t_{\ntoi} & \leq 1, \hspace{1.5cm} \forall \, n = 1,\ldots,N \label{eq:sr_t} \\
\hspace{-2cm} \lambda_n & \leq \mu_n, \, \mu_n  \leq 1, \,\forall \, n = 1,\ldots,N \\ 
\hspace{-2cm} P_n & \leq \Pnmax, \hspace{.9cm} \forall \, n = 1,\ldots,N \label{eq:sr_p} \\
\hspace{-2cm} \iRmax & \geq \sum_{n \in \cN_R} \mu_n e^{\kappa P_n + \kappa m_{n \rightarrow R} + \frac{\kappa^2}{2} \sigma^2_{n \rightarrow R}}, \nonumber \\
& \hspace{2.4cm} \forall \, R = 1,\ldots, N_R \label{eq:sr_interf}
\end{align}
%
%
\end{subequations}

\vspace{.3cm}

\noindent with $\{r_{n \rightarrow j}\}$ given by~\eqref{eq:rni}, and $\Pnmax := 10 \log_{10} \pnmax$. 
 
The non-convexity of constraints~\eqref{eq:sr_flow} and~\eqref{eq:sr_interf} makes problem (P1) non-convex, and thus hard to solve. Furthermore, function $Q(\cdot)$ in~\eqref{eq:rni} is difficult to handle in an optimization problem. In the next section, an approximate but efficiently solvable version of (P1) will be formulated. But first, some remarks are in order. 

\vspace{.3cm}

\noindent \textbf{Remark 1} \emph{(Monotonically non-decreasing utilities).} It follows from~\cite[Thm. 5]{Ribeiro}, that (P1) is optimally solved by setting $\{\lambda_n = \mu_n\}$ if each utility function $\cU_n(\rho_n)$ is monotonically non-decreasing. As many practical utilities satisfy this condition, $\cU_n(\rho_n)$ will be hereafter assumed non-decreasing, and variables $\{\lambda_n\}$ will be dropped. Strictly speaking, the choice $\{\lambda_n = \mu_n\}$ will lead to a solution of (P1) where queues may or may not be stable~\cite{Loynes}. On the other hand, condition $\lambda_n < \mu_n$ is challenging because it entails an open constraint set. From a practical perspective, queue stability can be readily ensured by imposing in~\eqref{eq:sr_p} the condition $\lambda_n + \epsilon \leq \mu_n$, with $0 < \epsilon \ll 1$ small enough, and replacing variables $\{\mu_n\}$ with $\{\bar \mu_n := \mu_n - \epsilon\}$. \hfill $\Box$

\noindent \textbf{Remark 2} \emph{(Conventional multi-hop networks).} The proposed routing framework can be considered also for non-CR multihop random access networks when node-to-node channels can not be estimated accurately - what could emerge with  e.g., a mobile ad hoc topology. Optimal routes and link reliabilities can be obtained by solving (P1), after discarding the interference constraints~\eqref{eq:sr_interf}, and re-defining the signal-to-noise ratio (SNR) of link $U_n \rightarrow U_i$ as $\gamma_{\ntoi} = p_n g_{\ntoi}/\sigma_i^2$. 

\hfill $\Box$

\noindent \textbf{Remark 3} \emph{(MAC protocol).} Since a random access protocol is adopted, a packet is as usual deemed lost when collisions among CR transmissions occur [cf.~\eqref{eq:rni}], and no mutual interference is explicitly modeled in~\eqref{eq:sinr}. However, the solution approach presented in the ensuing section can be effectively employed when different MAC strategies such as, e.g., carrier sensing medium access and orthogonal transmissions are utilized by the CR nodes; see also~\cite{Marques_book}. \hfill $\Box$

\section{Tractable Routing Protocol}
\label{sec:tractable}

To convexify constraint~\eqref{eq:sr_interf} it suffices to consider the
logarithmic change of variables $\tilde \mu_n := \ln(\mu_n)$. As for the flow constraint~\eqref{eq:sr_flow}, consider first introducing auxiliary variables $\{\nu_n\}$ representing the probability of CRs to remain silent, together with the extra constraints $\mu_n + \nu_n = 1$, for $n = 1,\ldots,N$. Further, a simple way to obtain a tractable approximation of $Q(x)$ consists in exploiting the commonly used upper and lower bounds proposed in~\cite{Chianitwc03,Haggman04}, which are very tight for $x > \sqrt{2}/2$.
Taking advantage of these bounds, and performing again a logarithmic change of variables $\tilde \nu_n := \ln(\nu_n)$, the probability $r_{\ntoi}$ can be (tightly) bounded as       
\begin{align}
r_{\ntoi} & \geq e^{\sum_{j \in \cI_{ni}} \tilde \nu_i} \times \left(1- \frac{1}{12}
e^{- \frac{1}{2} \left(\frac{P_n + m_{\ntoi} - \bar \Gamma_{\ntoi} }{\sigma_{\ntoi}} \right)^2 } \right. \nonumber \\
& \hspace{2.5cm}\left. - \frac{1}{4}
e^{- \frac{2}{3} \left(\frac{P_n + m_{\ntoi} - \bar \Gamma_{\ntoi} }{\sigma_{\ntoi}} \right)^2 }  \right)  \label{eq:lowerboundQ} \\
r_{\ntoi}  & \leq e^{\sum_{j \in \cI_{ni}} \tilde \nu_i} \times \left(1- \alpha_1
e^{- \alpha_2 \left(\frac{P_n + m_{\ntoi} - \bar \Gamma_{\ntoi} }{\sigma_{\ntoi}} \right)^2 } \right) \label{eq:upperboundQ}
\end{align}
where $\alpha_1 = 0.28$, and $\alpha_2 = 0.64$~\cite{Haggman04}. The premise for adopting the aforesaid bounds is that the decoding rate of CR links is at least $\approx 0.7$. This condition is met in practice if CRs and PUs are sufficiently far apart (see, e.g.~\cite{Shin04}). Furthermore, maximum packet error rates required for data and speech transmissions are generally considerably lower than $0.3$~\cite{Ephremides-twc02}. 

Consider now using the upper bound~\eqref{eq:upperboundQ} for the incoming traffic, and the lower bound~\eqref{eq:lowerboundQ} for the outgoing flows. As~\eqref{eq:lowerboundQ} and \eqref{eq:upperboundQ} are tight, this replacement not only  yields a tractable optimization problem, but also does not sacrifice optimality of the outcoming rates. With the logarithmic change of variable $\tilde t_{\ntoi} = \ln(t_{\ntoi})$, and after introducing auxiliary variables $\{\check y_{\ntoi} \geq \sqrt{2}/2\}$ and $\{\hat y_{\ntoi} \geq \sqrt{2}/2\}$, constraint~\eqref{eq:sr_flow} can be approximated as 
\begin{align}
\rho_n & + \sum_{i \in \cN_{n\rightarrow}} e^{\tilde \mu_n + \tilde t_{\ntoi} + \sum_{m \in \cI_{ni}} \tilde \nu_m } \left(  \frac{1}{12} e^{- \frac{1}{2} \hat y_{\ntoi}} +  \frac{1}{4} e^{- \frac{2}{3} \hat y_{\ntoi}}   \right)  \nonumber \\ 
& \hspace{-.5cm}  +  \sum_{j \in \cN_{\rightarrow n}} e^{\tilde \mu_j + \tilde t_{\jton}  + \sum_{m \in \cI_{jn}} \tilde \nu_m }  - \sum_{i \in \cN_{n \rightarrow}} e^{\tilde \mu_n + \tilde t_{\ntoi} + \sum_{m \in \cI_{ni}} \tilde \nu_m }  \nonumber \\
& \hspace{-.5cm}  - \alpha_1  \sum_{j \in \cN_{\rightarrow n}}  e^{\tilde \mu_j + \tilde t_{\jton} + \sum_{m \in \cI_{jn}} \tilde \nu_m - \alpha_2 \check y_{\jton}} \leq 0 
\label{eq:approx_flow} 
\end{align}
with the auxiliary constraints
\begin{align}
\sigma_{\ntoi}^{\frac{1}{2}} (\hat y_{\ntoi}) \leq P_n + m_{\ntoi} - \bar \Gamma_{\ntoi} \label{eq:constraint_yout} \\ 
\sigma_{\iton}^{\frac{1}{2}} (\check y_{\iton}) \geq P_i + m_{\iton} - \bar \Gamma_{\iton} \label{eq:constraint_yin} .
\end{align}

For notational convenience, define the variable vector $\bx_{n} := [P_n, \{\tilde t_{\ntoi}\},\tilde \mu_n, \tilde \nu_n, \{\check y_{\jton}, \hat y_{\ntoi}\}]^T$ per node $U_n$, for $n = 1,\ldots,N$.  Upon re-expressing in a compact form the flow constraint~\eqref{eq:approx_flow} as $f_n(\{\bx_n\}) \leq 0$, and defining the constraint set $\cB_{\bx_n}$ per node $U_n$ as
\begin{align}
\cB_{\bx_n} & := \left\{ \bx_n: \sum_{i \in \cN_{n\rightarrow}} e^{\tilde t_{\ntoi}} \leq 1, e^{\tilde \mu_n} + e^{\tilde \nu_n} \leq 1, \right. \nonumber \\ 
& \hspace{2.5cm} P_n\leq \Pnmax, \textrm{and } \eqref{eq:constraint_yout}, \eqref{eq:constraint_yin} \Big\}
\end{align}
where (non)negativity of the variables is left implicit, problem (P1) can be re-formulated as  
\begin{subequations}
%
%
\begin{align}
 \textrm{(P2)} \quad  \max_{\{\bx_n \in \cB_{\bx_n}\}} 
 \,\,\, & \sum_{n=1}^N \cU_n(\rho_n) - \sum_{n=1}^N \cC_n(P_n) \label{eq:sr_cost2} \\
\textrm{subject to } f_n(\{\bx_n\}) & \leq 0,  \hspace{.6cm} n=1,\ldots, N \label{eq:flow_compact} \\
\sum_{n \in \cN_R} f_{\iota_R,n}(\bx_n) & \leq \iRmax, \, R=1,\ldots, N_R \label{eq:sr_interf2} 
\end{align}
%
%
\end{subequations}
where $f_{\iota_R,n}(\bx_n) := e^{\tilde \mu_n + \kappa P_n + \kappa m_{n \rightarrow R} + \frac{\kappa^2}{2} \sigma^2_{n \rightarrow R}}$.

Constraints~\eqref{eq:approx_flow} are still non-convex because the last two sums (with their signs) are concave, and likewise~\eqref{eq:constraint_yout} is also concave. Nevertheless, the structure of (P2) allows convex approximation methods for obtaining its solution efficiently. Among candidate methods, the successive convex approximation approach~\cite{MaWr78} is well suited for the problem at hand because it guarantees first-order KKT optimality under mild regularity conditions.

\subsection{KKT solution via successive convex approximation}
\label{sec:sca}

The general successive convex approximation method is outlined first. Suppose that the objective function to be maximized is concave in the optimization variables $\bx$, and the constraint set is the intersection of a set $\cA := \{\bx| f_n(\bx) \leq 0, n = 1,\ldots,N \}$ with a convex set $\cB$, which captures convex constraints, if any. Assume that $f_n(\bx)$, $n=1,\ldots,N$, are differentiable but generally non-convex functions. Then, starting from a feasible point $\bx^{(0)} \in \cA \cap \cB$, a
series $\ell = 1,\ldots$, of surrogate problems is solved, where $\cA$ is substituted per iteration $\ell$ by a convex set $\cA^{(\ell)}$. Since the intersection of convex sets yields a convex set, the resulting optimization problems are convex. 
For each $n = 1,\ldots,N$, let $\tilde f_n(\bx; \bx^{(\ell)})$ denote the surrogate convex function for $f_n(\bx)$, which may depend on the solution $\bx^{(\ell)}$ to the problem of the previous $(\ell-1)$-st iteration. Then, the convex set $\cA^{(\ell)}$ is constructed as $\cA^{(\ell)} := \{\bx| \tilde f_n(\bx;\bx^{(\ell)}) \leq 0, n = 1,\ldots,N \}$. 
Provided that each function $\tilde f_n(\bx;\bx^{(\ell)})$, $n = 1,\ldots,N$, is convex, differentiable, and satisfies conditions~\cite{MaWr78} 
\begin{enumerate}
\item[c1)] $f_n(\bx) \leq \tilde f_n(\bx;\bx^{(\ell)}),\quad \forall \bx \in \cA^{(\ell)} \cap \cB$

\item[c2)] $f_n(\bx^{(\ell)}) = \tilde f_n(\bx^{(\ell)};\bx^{(\ell)})$, and

\item[c3)] $\nabla f_n(\bx^{(\ell)}) = \nabla \tilde f_n(\bx^{(\ell)};\bx^{(\ell)})$ 
\end{enumerate}
the series of solutions to the approximate problems converge to the KKT point of (P2)~\cite{MaWr78}.

In order to apply the successive convex approximation method to
(P2), surrogate constraints for the non-convex constraints must be determined. The first three terms in~\eqref{eq:approx_flow} are convex, whereas the fourth and fifth terms are concave. Letting $- e^{x_1 + \beta x_2 - \alpha x_3}$ represent one of the non-convex summands, a convex surrogate function satisfying c1)-c3) can be obtained by replacing the non-convex summands with the affine function 
\begin{align}
&\hspace{-.4cm} - e^{x_1 + \beta x_2 - \alpha x_3}  
\leq e^{x_1^{(\ell)} + \beta x_2^{(\ell)} - \alpha x_3^{(\ell)}}  \times  \left[(x_1^{(\ell)} - x_1) \right. \nonumber \\ 
&  \hspace{2.3cm}  \left. + \beta (x_2^{(\ell)} - x_2) - \alpha (x_3^{(\ell)} - x_3) -1 \right] \,.
\label{eq:upper_exp} 
\end{align}
As for~\eqref{eq:constraint_yout}, an
upper-bound of $\sqrt{\hat y_{\ntoj}}$ can be obtained via the supporting hyperplane, and the resulting surrogate convex constraints become
\begin{align}
\frac{\hat y_{\ntoj} - \hat y_{\ntoj}^{(\ell)}}{2 \sqrt{\hat y_{\ntoj}^{(\ell)}}} + \sqrt{\hat y_{\ntoj}^{(\ell)}} - P_n - m_{\ntoi} + \bar \Gamma_{\ntoi} \leq 0 .
\label{eq:upper_sqrt} 
\end{align}

Overall, the problem to solve in the $\ell$-th iteration is
given by (P2) with~\eqref{eq:constraint_yout}
replaced by~\eqref{eq:upper_sqrt} to form the surrogate constraint set $\tilde \cB_n$, and by employing~\eqref{eq:upper_exp} along the feasible points $\{\bx_n^{(\ell)}\}_{n=1}^{N}$ 
to obtain a surrogate convex flow conservation constraint $\tilde f_n^{(\ell)}(\{\bx_n\}) \leq 0$; that is, 
\begin{subequations}
\begin{align}
 (\textrm{P2}^{(\ell)}) \quad\quad  \max_{ \{\bx_n \in \tilde \cB_n\}} 
 \,\,\, \sum_{n=1}^N \cU_n(\rho_n) - & \sum_{n=1}^N \cC_n(P_n) \label{eq:sr_cost3} \\
\textrm{subject to } \eqref{eq:sr_interf2}~\textrm{and}~ 
\tilde f_n^{(\ell)}(\{\bx_n\}) & \leq 0, \, \forall \, n \label{eq:flow_compact_ell} 
\end{align}
\end{subequations}  
Problem $(\textrm{P2}^{(\ell)})$ is convex, and thus efficiently solvable using interior-point methods~\cite{BoV04}. It is worth mentioning that the solution of (P2$^{(\ell)}$), $\ell = 1,2,\ldots$ always lies inside the feasibility region of the original non-convex problem (P2)~\cite{MaWr78}. This observation suggests readily an \emph{online} practical implementation of the algorithm whereby node $U_n$ does not necessarily wait for the successive convex approximation algorithm to converge, but rather relies on $\bx_n^{(\ell)}$ as and when it becomes available. In the limit (i.e., for $\ell \gg 1$), $\bx_n^{(\ell)}$ will be KKT-optimal. An online implementation of the iterative optimization allows tracking of slow variations in the network topology and SINR statistics.

\section{Distributed statistical routing}
\label{sec:distributed}

To obviate the high communication cost associated with the collection of channel statistics for all links at a central processing unit, and the subsequent dissemination of the optimized variables, it is of prime interest to solve (P2) in a distributed manner. A distributed cross-layer optimization algorithm is also desirable because of its scalability with regards to power requirements and network size, and robustness to isolated points of failure. 

Distributing (P2) is tantamount to developing a distributed solver for each of the convex problems (P2$^{(\ell)}$), $\ell = 1,2,\ldots$. To this end, it is necessary to decompose (P2$^{(\ell)}$) into smaller sub-problems, which can be locally solved by nodes $\{U_n\}$ via \emph{local} message exchanges. Unfortunately, the interference constraints~\eqref{eq:sr_interf2} challenge decomposability,
as they couple portions of the CR network. Furthermore, for each $U_n$, constraint~\eqref{eq:flow_compact_ell} involves variables pertaining to the one-hop neighboring nodes $U_i \in \cN_{\rightarrow n}$, and to CRs in the collision-related sets  $\{\cI_{jn}\}$ and $\{\cI_{ni}\}$. To overcome the first hurdle, consider first the following problem 
\begin{subequations}
\begin{align}
 (\textrm{P3}^{(\ell)})\quad g(\{\inRmax(\ell)\}) := \nonumber \\
 & \hspace{-2cm}  \max_{\{\bx_n \in \tilde \cB_n \}} 
 \,\,\,  \sum_{n=1}^N \cU_n(\rho_n) - \sum_{n=1}^N \cC_n(P_n) \label{eq:sr_costP3} \\
\textrm{subject to }
\tilde f_n(\{\bx_n\}) & \leq 0,  \hspace{1cm} \forall \, n \label{eq:flow_dist} \\
f_{\iota_R,n}(\bx_n) & \leq \inRmax(\ell), \,\, \forall \, n, R 
 \label{eq:interf_dist3} 
\end{align}
\end{subequations}
where the interference $\iRmax$ for PU $R$ is pre-partitioned in \emph{given} per-CR fractions $\{\{\inRmax(\ell)\}\}_{n \in \cN_R}$. Problem $(\textrm{P2}^{(\ell)})$ will be revisited later on. Then, collect \emph{local} copies of $\bx_{j \rightarrow n} := [\tilde t_{\jton}, P_j, \tilde \mu_j]^T$ at node $U_n$, for each $j \in \cN_{\rightarrow n}$ into a vector $\bx_{j \rightarrow n,n} := [\tilde t_{\jton,n}, P_{j,n}, \tilde \mu_{j,n}]^T$. 
Likewise, let $\{\tilde \nu_{n,m}\}$ denote local copies of $\{\tilde \nu_m | m \in \cI_n  \}$, with $\cI_n := (\cup_{i \in \cN_{n \rightarrow}} \cI_{ni}) \cup ( \cup_{j \in \cN_{\rightarrow n}} \cI_{jn})$; i.e., local copies of $\tilde \nu_m$ for users that may interfere with $U_n$'s transmissions. Then, $(\textrm{P3}^{(\ell)})$ can be equivalently re-formulated as  
\begin{subequations}
\begin{align}
 (\textrm{P4}^{(\ell)}) \quad\quad  \max_{\{\bx_n \in \tilde  \cB_{\bx_n} \}} 
 \,\,\, \sum_{n=1}^N \cU_n(\rho_n) - & \sum_{n=1}^N \cC_n(P_n) \label{eq:sr_costP4} \\
 & \hspace{-4cm} \textrm{subject to } \nonumber \\
\tilde f_n(\bx_n,\{\bx_{j\rightarrow n,n}\},\{\nu_{m,n}\}) & \leq \,\,  \forall \, n\label{eq:flow_dist} \\
& \hspace{-4cm} f_{\iota_R,n}(\bx_n) \leq \inRmax(\ell), \quad  \forall \, n, R 
 \label{eq:interf_dist4} \\ 
& \hspace{-3.4cm} \bx_{j\rightarrow n}  = \bx_{j\rightarrow n,n},  \quad j \in \cN_{\rightarrow n},   \forall \, n, \label{eq:equalities_x} \\
& \hspace{-3.0cm} \tilde \nu_{m}  = \tilde \nu_{m,n}, \quad \,\,\, \, m \in \cI_{n},  \forall \, n  \label{eq:equalities_nu} 
\end{align}
\end{subequations}
where the notation $\tilde f_n(\bx_n,\{\bx_{j\rightarrow n,n}\},\{\nu_{m,n}\})$ emphasizes the dependence of the surrogate flow conservation constraint $\tilde f_n(\cdot)$ on the newly introduced local variables. Problem $(\textrm{P4}^{(\ell)})$ is amenable to a distributed solution, where~\eqref{eq:equalities_x}-\eqref{eq:equalities_nu} can be enforced by means of local message passing.  

Suppose that there is a non-zero probability (possibly multi-hop) directed path connecting $U_n$ to  nodes $U_m \in \cN_{\rightarrow n} \cup \cN_{n \rightarrow} \cup \cI_n$; i.e., nodes coupled in the optimization problem. If not, a control channel can be employed as usual. Problem $(\textrm{P4}^{(\ell)})$ may be solved in a distributed manner using the dual sub-gradient method~\cite{Ribeiro-twc08,Palomar06}. However, recovering the primal variables $\{\bx_n\}$ from the Lagrange multipliers optimizing the dual function is not always guaranteed if the objective in~\eqref{eq:sr_costP4} is not strictly convex, and the step-size in the sub-gradient ascent is constant. Furthermore, primal averaging can not be performed in this case, unless the equality constraints are appropriately relaxed~\cite{Nedic09}.

\begin{algorithm}[t]
\label{alg:ADMoM}
\caption{Distributed algorithm for (P4$^{(\ell)}$)} \small{
\begin{algorithmic}
\STATE Assumption: bidirectional links, or bidirectional control channels.
\STATE Use solution of (P4$^{(\ell-1)}$) to initialize variables. 

\FOR {$l=0,1,\ldots$ (repeat until convergence)} 

\STATE Receive multipliers $\{\bq_{n,i}(l)\}_{i \in \cN_{n \rightarrow}}$, $\{v_{n,p}(l)\}_{p \in \{p | U_n \in \cI_{p}\}}$.
\STATE Update $\bar \bx_n := \{ \bx_n, \{\bx_{j \rightarrow n,n}\}, \{\nu_{m,n}\}\}$ via~\eqref{eq:iterI1}  
\STATE Transmit $\bx_{n \rightarrow i}(l)$ to $U_i \in \cN_{n \rightarrow}$, and $\bx_{j \rightarrow n,n}(l)$ to $U_j \in \cN_{\rightarrow n}$

\STATE Transmit $\{\tilde \nu_{m,n}\}$ to $U_m \in \cI_n$ via neighboring nodes. Forward $\{\tilde \nu_{m,j}\}_{j \in \cN_{\rightarrow n}}$ to $U_i \in \cN_{n \rightarrow}$.

\STATE Receive $\bx_{j \rightarrow n}(l)$ from $U_i \in \cN_{n \rightarrow}$, and $\bx_{n \rightarrow i,i}(l)$ from $U_i \in \cN_{n \rightarrow}$

\STATE Receive $\tilde \nu_{n,p}(l)$ from $U_p \in \{U_p | U_n \in \cI_{p}\}$

\STATE Dual update via~\eqref{eq:update_multipliers_q}-\eqref{eq:update_multipliers_v}. 

\STATE Transmit multipliers $\bq_{n,i}(l+1)$ to $U_i \in \cN_{n \rightarrow}$.

\STATE Transmit $v_{m,n}(l+1)$ to $U_m \in \cI_n$ via neighboring nodes, forward $\{v_{m,j}\}_{j \in \cN_{\rightarrow n}}(l+1)$ to $U_i \in \cN_{n \rightarrow}$.  

\STATE Use parameters $\bar \bx_n$ to transmit data in case of on-line implementation.  

\ENDFOR
\end{algorithmic}}
\end{algorithm}

One effective remedy is offered by the alternating direction method of multipliers (ADMoM), where the optimization argument in $(\textrm{P4}^{(\ell)})$ is augmented with a quadratic regularization term corresponding to the squared norm of the equality constraints~\cite[Sec.~3.4]{BeT89}. Specifically, letting $\{\bq_{j,n}\}$ and $\{v_{m,n}\}$ denote the multipliers associated with the equality constraints~\eqref{eq:equalities_x} and \eqref{eq:equalities_nu}, respectively, the partial quadratically-augmented Lagrangian function is given by
\begin{align}
& \cL( \{ \bar \bx_n\}, \{\bq_{j,n}\},\{v_{m,n}\}):=  \sum_{n=1}^{N } \Big[ -\cU_n(\rho_n) +  \cC_n(P_n)  \nonumber \\ 
 &   + \sum_{j \in \cN_{\rightarrow n}} \left(\bq_{j,n}^T(\bx_{j \rightarrow n} - \bx_{j \rightarrow n,n})  + \frac{c}{2} \| \bx_{j \rightarrow n} - \bx_{j \rightarrow n,n}\|_2^2 \right)  \nonumber \\ 
 & + \sum_{m \in \cI_{n}} \left( v_{m,n}( \tilde \nu_{m} - \tilde \nu_{m,n} )  + \frac{c}{2} ( \tilde \nu_{m} - \tilde \nu_{m,n} )^2  \right) \Big]
\label{eq:lagrangian}  
\end{align}  
where $\bar \bx_n := \{ \bx_n, \{\bx_{j \rightarrow n,n}\}, \{\nu_{m,n}\}\}$, and $c>0$ is an arbitrary constant. Notice that $\cL(\cdot)$ is defined over the primal feasible region $\cA := \cap_{n=1}^{N} \cA_{n}$, with $\cA_{n} := \{ \bar \bx_n  | \bx_n \in \tilde  \cB_{\bx_r},~\eqref{eq:flow_dist}, \eqref{eq:interf_dist4} \}$. ADMoM amounts to performing the following iterations ($l$ denotes the iteration index) 

\vspace{.2cm}

\noindent \textbf{[I.1]} \emph{Primal update}. Given $\{\bq_{j,n}(l)\}$ and $\{v_{m,n}(l)\}$, update primal variables in a coordinate descent fashion; i.e., for $n=1,\ldots,N$, update  $\bar \bx_n$ as:
\begin{subequations}
\begin{align}
& \bar  \bx_n(l+1) := \min_{ \bar  \bx_n \in \cA_{n} } 
 \,\,\, 
 \cL_n(\bar \bx_n,l) 
 \label{eq:iterI1} \\ 
& \cL_n(\bar \bx_n,l) := \cL( \bar \bx_1(l),\ldots, \bar \bx_{n-1}(l), \bar \bx_n, \bar \bx_{n+1}(l),\ldots, \bar \bx_N(l), \nonumber  \\
& \hspace{2cm}  \{\bq_{j,n}(l)\},\{v_{m,n}(l)\})
\end{align}
where $\cL_n(\bar \bx_n,l)$ is obtained by keeping $\{\bar \bx_m(l) \}_{m \neq n}$ fixed to their values at iteration $l$. 

\vspace{.2cm}

\noindent \textbf{[I.2]} \emph{Dual update}. Given the primal variables $\{\bar \bx_n(l+1)\}_{n = 1}^{N}$, updated multipliers as:
\begin{align}
\bq_{j,n}(l+1) & = \bq_{j,n}(l) + \beta \left[\bx_{j \rightarrow n}(l+1) - \bx_{j \rightarrow n,n}(l+1) \right]  \label{eq:update_multipliers_q}   \\
v_{m,n}(l+1) & = v_{m,n}(l) + \beta\left[\tilde \nu_{m}(l+1) - \tilde \nu_{m,n}(l+1) \right] \label{eq:update_multipliers_v}   
\end{align}
\end{subequations}
where $\beta > 0$ is the step-size.

\vspace{.2cm}

Once the primal iterates of the neighboring nodes $\{\bx_{j\rightarrow n}(l+1)\}$ and $\{\tilde \nu_{m}(l+1)\}$ become available at node $U_n$, the dual updates~\eqref{eq:update_multipliers_q}-\eqref{eq:update_multipliers_v} can be performed locally. As for the primal update, the local augmented Lagrangian  [cf.~\eqref{eq:lagrangian}]
\begin{align}
& \cL_n(\bar \bx_n,l) = -\cU_n(\rho_n) +  \cC_n(P_n)  \nonumber \\ 
& + \sum_{j \in \cN_{\rightarrow} n} \left[-  \bq_{j,n}^T(l)\bx_{j \rightarrow n,n} + \frac{c}{2} \| \bx_{j \rightarrow n}(l) - \bx_{j \rightarrow n,n}\|_2^2 \right] \nonumber \\ 
& + \sum_{i \in \cN_{n \rightarrow}} \left[ \bq_{n,i}^T(l) \bx_{n \rightarrow i}  + \frac{c}{2} \| \bx_{n \rightarrow i} - \bx_{n \rightarrow i,i}(l)\|_2^2  \right] \nonumber \\
& + \sum_{m \in \cI_{n}} \left[ - v_{m,n}(l) \tilde \nu_{m,n} + \frac{c}{2} \left( \tilde \nu_{m}(l) - \tilde \nu_{m,n} \right)^2  \right]
 \nonumber \\  
& + \sum_{p | n \in \cI_{p}} \left[ v_{n,p}(l) \tilde \nu_{n} + \frac{c}{2} \left( \tilde \nu_{n} - \tilde \nu_{n,p}(l) \right)^2   \right]
\label{eq:Lagrangian_local} 
\end{align}
can be minimized at node $U_n$ upon collecting $\bx_{n \rightarrow i,i}(l)$ and multipliers $\{\bq_{n,i}(l)\}$ from the one-hop neighboring nodes $U_i \in \cN_{n \rightarrow}$, and $\{v_{n,p}(l)\}$ and $\tilde \nu_{n,p}(l)$ from nodes $U_p \in \{U_p | U_n \in \cI_{p}\}$; that is, from the nodes whose transmissions can collide with the ones of $U_n$. Roughly speaking, the latter quantities pertain to the two-hop neighborhood of node $U_n$ and are due to the basic properties of the random access strategy. If a different medium access protocol such as, e.g., CSMA is employed, \eqref{eq:equalities_nu} will not be required and the message-passing overhead can be further reduced.  The ADMoM-based distributed algorithm is tabulated as in Algorithm~1, and the convergence to the optimal primal arguments $\{\bar \bx_n(l)\}$  as $l \rightarrow \infty$ is summarized in the following proposition.

\begin{proposition}
If there exists a non-zero probability (possibly multi-hop) directed path connecting $U_n$ to nodes $U_m \in \cN_{\rightarrow n} \cup \cN_{n \rightarrow} \cup \cI_n$, for all $n$, the iterates $\{\bar \bx_n(l)\}$ generated by Algorithm~1 converge to a globally optimal solution to $(\textrm{P4}^{(\ell)})$.  
\end{proposition}

\emph{Proof}. Existence of a path connecting $U_n$ to nodes $U_m \in \cN_{\rightarrow n} \cup \cN_{n \rightarrow} \cup \cI_n$, guarantees a regular exchange of local primal variables and multipliers among neighboring nodes. Under this assumption, convergence of the primal iterates $\{\bar \bx_n(l)\}$ to their optimal values as $l \rightarrow \infty$ can be readily established using the result in~\cite[Prop.~4.2]{BeT89}. \hfill $\Box$   

Algorithm~1 can also be implemented in an online fashion. The equality constraint violation during the initial iterations of the algorithm may induce an initial increase of some queues. Thus, an online implementation is feasible if nodes can afford such a potential increase in the queue length before reaching consensus on the local variables.

\subsection{Handling the interference constraint via primal decomposition}
\label{sec:primal}

Reconsider now problem (P2$^{(\ell)}$), where the interference budgets $\{\iRmax\}_{R = 1}^{N_R}$ are \emph{not} partitioned \emph{a priori} among CR nodes. As primal variables become feasible only when dual decomposition algorithms have converged,  utilization of network parameters obtained from intermediate iterates can possibly lead to violation of the interference constraint. 
To enforce strict PU protection during network operation, the primal decomposition technique is invoked here; see, e.g.,~\cite{Palomar06}. With this method, resources shared among CR nodes are essentially allocated by a master problem. Specifically, at each iteration $k = 1, 2,\ldots$ of the primal decomposition algorithm, problem (P3$^{(\ell,k)}$) is solved for given $\{\{\inRmax(\ell,k)\}_{R,n}\}$; then, $\{\{\inRmax(\ell,k)\}_{R,n}\}$ are updated by solving the following problem:
\begin{subequations}
\begin{align}
 (\textrm{P5}^{(\ell,k)}) \quad\quad \{\inRmax(\ell,k)\} =  \arg  \max_{\{\inRmax\}} & g(\{\inRmax\}) \\
\textrm{subject to } \inRmax & \geq 0, \\
\quad \sum_{n \in \cN_R} \inRmax & \leq \iRmax, \, \forall \,\, R \label{eq:sum_interf} .
\end{align}
\end{subequations}
To solve $(\textrm{P5}^{(\ell,k)})$, the subgradient algorithm can be employed~\cite{BoV04}. Specifically, the subgradient of $g(\{\inRmax(\ell,k)\})$ with respect to $\inRmax(\ell,k)$ is given by the \emph{optimal} Lagrange multiplier $u_{R,n}(k)$ corresponding to the constraint $f_{\iota_R,n}(\bx_n) \leq \inRmax(\ell,k)$ in (P3$^{(\ell)}$) at iteration $k$~\cite{Palomar06}. Therefore, $\inRmax(\ell,k)$ is updated as
\begin{align}
\inRmax(\ell,k+1) = \cP_{\iota_R} \left\{\inRmax(\ell,k) + \xi(k+1) u_{R,n}(k)  \right\} 
\label{eq:projection} 
\end{align}
where $\xi(\cdot)$ is the step size, and $\cP_{\iota_R}\{\cdot\}$ denotes projection onto the region defined by~\eqref{eq:sum_interf}, operation that can be efficiently computed as in, e.g.~\cite{Michelot86}. At each step $k$ of the primal algorithm, CR nodes can employ variables obtained from (P3$^{(\ell,k)}$) for network operation, as PU interference protection is enforced by updates~\eqref{eq:projection}. 

The projection in~\eqref{eq:projection} needs to be performed by a ``head node'' in the CR sub-network CR $\cN_R$, which is formed by nodes that are coupled by constraint $\sum_{n \in \cN_R} f_{\iota_R,n}(\bx_n) \leq \inRmax(\ell)$. Per iteration $k$, the head node has to collect the optimal Lagrange multipliers from the CRs in $\cN_R$, and then broadcast the updated interference budgets $\{\inRmax(\ell,k+1)\}_{n}$. This leads to a semi-distributed algorithm, but the high message-passing overhead entailed by centralized solutions is nonetheless alleviated.  
The online algorithm obtained through the successive convex approximation and primal decomposition is tabulated in Algorithm~2, and its convergence properties are summarized next.      

\begin{proposition}
If there is a cycle connecting nodes $U_n \in \cN_R$, $\forall \,\, R$, the iterates generated by Algorithm 2 converge to a KKT solution to $(\textrm{P2})$.  
\end{proposition}

\emph{Proof.} Since the original problem (P2$^{(\ell,k)}$) is convex, the subproblems (P3$^{(\ell,k)}$) as well as the master problem (P5$^{(\ell,k)}$) are all convex, and thus the globally optimal solution of (P2$^{(\ell)})$ is attained via primal decomposition~\cite{Palomar06}. Existence of a cycle connecting nodes $U_n \in \cN_R$ ensures that the multipliers $\{\xi(k)\}$ can be collected to a cluster head node, and that $\{\inRmax(\ell,k+1)\}$ can be subsequently sent back. Finally, since (P2$^{(\ell)})$ is optimally solved per iteration $\ell$ of the successive convex approximation, convergence of Algorithm~2 to a KKT point of (P2) is guaranteed~\cite{MaWr78}. \hfill $\Box$    \\

\noindent \textbf{Remark 4} \emph{(Fully distributed algorithm)}. At the expense of possibly sacrificing optimality of the resultant exogenous rates, coefficients $\{\inRmax\}$ can be set a priori based on the distance between CRs and PU $R$. This may be reasonable especially if shadowing can not be estimated~\cite{DKGtvt11}. In this case, it is not necessary to compute the primal decomposition iterates. \hfill $\Box$

\begin{algorithm}[t]
\label{alg:primal}
\caption{Overall on-line algorithm for (P2)} \small{
\begin{algorithmic}
\STATE Assumption: Path connecting all $U_n \in \cN_R$, for all PUs $R = 1, \ldots, N_R$.

\FOR {$\ell=1,\ldots$ (repeat until convergence)} 

\STATE Use solution of (P2$^{(\ell-1)}$) to compute~\eqref{eq:upper_exp}-\eqref{eq:upper_sqrt}. If $\ell = 1$, use suitable feasible point.

\FOR {$k=0,1,\ldots$ (repeat until convergence)} 

\STATE Receive $u_{R,n}(k)$ from head node.

\STATE Solve (P3$^{(\ell,k)}$) using Algorithm~1.  

\STATE Transmit multiplier $u_{R,n}(k)$ to head node via neighboring nodes.

\STATE If head node: update $\{\inRmax(\ell,k+1)\}$ via~\eqref{eq:projection}. 

\STATE Utilize $\bar \bx_n(k)$ for network operation.

\ENDFOR

\ENDFOR
\end{algorithmic}}
\end{algorithm}

\section{Packet deliverability in dynamic CR environments}
\label{sec:deliverability}

Statistics of the SINR may vary during network operation, because of the dynamic nature of shadow fading~\cite{DKGtvt11}, and the variable PU interference levels [cf.~\eqref{eq:sinr}]. CR topology may also change with time. Proximity of PUs with intermittent activity, or, mobile PU devices may loose link connectivity during certain time intervals. The routing problem (P2) must be (re-)solved whenever network topology and SINR statistics change. Alternatively, it can be implemented online to track slow environmental dynamics. Either way, it is necessary to establish conditions ensuring that packets are eventually delivered to the sink when routes, MAC, and physical layer parameters are regularly updated. 

Let $s_n(\tau) \in \{0,1\}$ be a binary variable taking value $1$ if a packet, after having been randomly routed through the network, is placed in $U_n$'s queue at time $\tau$, and let $\vartheta_n(\tau) := \Prob\{s_n(\tau) = 1\}$ denote the probability of such an event. Further, collect $\{\vartheta_n(\tau)\}$ in the $(N + 1) \times 1$ vector $\bvartheta(\tau) := [\vartheta_1(\tau),\ldots,\vartheta_{N+1}(\tau)]^T$. 
CR-PU hierarchy may prevent CR nodes from forwarding packets during certain time intervals. Let $\ell_{\ntoj}$ be a binary variable that takes value $1$ if link $U_n \rightarrow U_j$ is active, and define $\chi_{\ntoj} := \Prob\{\ell_{\ntoj} = 1\}$. If active, link $U_n \rightarrow U_j$ is characterized by a link reliability $r_{\ntoj}(\tau)$. Probabilities $\{\chi_{\ntoj}\}$ clearly depend on PU activity factors and locations, and determine the average connectivity of the CR network~\cite{abbagnale-Secon}. 

If a packet is in $U_n$'s queue at time $\tau$, then $U_n$ may decide with probability $t_{\ntoj}(\tau)$ to route it through one of the available links, where index $\tau$ emphasizes the time-variability of routes. Clearly, if neither node locations nor the PU interference or channel conditions change for a certain number of time slots, then $\{t_{\ntoi}(\tau)\}$ and $\{r_{\ntoi}(\tau)\}$ remain invariant. The evolution of $\{\vartheta_{n}(\tau)\}$ can thus be fully characterized by the product probabilities $\{t_{\ntoj}(\tau) \, r_{\ntoj}(\tau)\}$, and the link availability factors $\{\chi_{\ntoj}\}$. Upon invoking the law of total probability, it holds that $\vartheta_{n}(\tau+1) =  \sum_{i = 1}^{N+1} \Prob\{s_n(\tau+1)=1|s_i(\tau)=1,\ell_{\iton} = 1\} \Prob\{s_i(\tau)=1\} \Prob\{\ell_{\iton} =1\}  = \sum_{i = 1}^{N+1} t_{\ntoj}(\tau) r_{\ntoj}(\tau) \chi_{\ntoj} \vartheta_i(\tau)$. 
Define the $(N+1) \times (N+1)$ packet delivery probability matrix $\bD(\tau)$, whose off-diagonal entry $(i,n)$ is $\{t_{n \rightarrow i} r_{n \rightarrow i} \chi_{n \rightarrow i}\}$ if $U_i$ is a one-hop neighbor of $U_n$, and $0$ otherwise. The diagonal entry $(n,n)$ of $\bD(\tau)$ represents the probability that a packet remains in $U_n$'s queue, which equals $1-\sum_{i\neq n} t_{n \rightarrow i} r_{n \rightarrow i} \chi_{n \rightarrow i}$. Finally, since the sink node will not route packets to any other node, set the $(n,N+1)$-th entry of $\bD(\tau)$ to $D_{n,N+1}(\tau) = 0$, and $D_{N+1,N+1}(\tau) = 1$. Matrix $\bD(\tau)$ is by construction a column stochastic, meaning that $\bD^T(\tau) \mathbf{1}_{N+1} = \mathbf{1}_{N+1}$ for all $\tau$. Then, the evolution of $\{\vartheta_{n}(\tau)\}$ can be expressed in matrix-vector form as $\bvartheta(\tau+1) = \bD(\tau) \bvartheta(\tau)$. 

Using an inductive argument, it is possible to show that the $(i,n)$th entry of the stochastic matrix $\bar \bD(t) := \prod_{\tau = 1}^{t} \bD(\tau)$ represents the probability that a packet generated at $U_n$ reaches node $U_i$ in $t$ time slots~\cite[Ch.~2]{Newman}. Therefore, it readily follows that a packet is eventually delivered to the sink node $U_{N+1}$ if and only if 
\begin{equation}
\lim_{t \rightarrow + \infty} \bvartheta(t) = \lim_{t \rightarrow + \infty} \bar \bD(t) \bvartheta(0) = [\mathbf{0}_{N}^T \, 1]^T
\label{eq:deliverability}
\end{equation}
holds for any initial distribution $\bvartheta(0)$. A simple condition on the CR network topology is provided next in order for~\eqref{eq:deliverability} to be satisfied.

\begin{proposition} If $\sum_{i \neq n} \chi_{\iton} > 0$, $\forall \,\, \{U_n\}_{n = 1}^{N}$ and $\sum_{n} \chi_{n \rightarrow N+1} > 0$, a packet stochastically routed according to probabilities $\{t_{n \rightarrow i}(\tau)\}$ over links with reliabilities $\{r_{n \rightarrow i}(\tau)\}$ will be eventually delivered to the destination with probability (w.p.) $1$.
\label{prop:deliverability}  
\end{proposition}
\noindent \emph{Proof.} The conditions of Prop.~\ref{prop:deliverability} ensure that there exists a multi-hop path connecting each node to the destination $U_{N+1}$ in the average connectivity graph, where link $U_n \rightarrow U_j$ is present if $\chi_{\ntoj} > 0$~\cite{abbagnale-Secon}. Let $t^*$ be the minimum number of time slots such that $U_{N+1}$ can be reached from any node with non-zero probability; i.e., $t^* = \min\{t: \bar D_{N+1,n}(t) > 0 \,\, \forall n = 1,\ldots, N\}$. Then, the probability that a packet is in $U_{N+1}$'s queue at time $t^* + 1$ is given by
\begin{equation}
\vartheta_{N+1}(t^*+1) = \sum_{n = 1}^N \bar D_{N+1,n}(t^*) \vartheta_n(t^*) + \vartheta_{N+1}(t^*) \, .
\label{eq:pdelivery}
\end{equation}
Arguing by contradiction, suppose that $\lim_{t^* \rightarrow + \infty} \vartheta_{N+1}(t^*) = \alpha < 1$; meaning that the packet is not delivered to $U_{N+1}$ w.p. $1-\alpha > 0$. Taking the limit on both sides of~\eqref{eq:pdelivery}, one arrives at 
\begin{align}
& \hspace{-.8cm} \lim_{t^* \rightarrow + \infty} \vartheta_{N+1}(t^*+1) \nonumber \\ 
& = \lim_{t^* \rightarrow + \infty} \left[ \sum_{n = 1}^N \bar D_{N+1,n}(t^*) \vartheta_n(t^*) + \vartheta_{N+1}(t^*) \right] \nonumber \\
& \geq \min_n \{ \bar D_{N+1,n}(t^*) \} \sum_{n = 1}^N \lim_{t^* \rightarrow + \infty} \vartheta_n(t^*)  + \alpha .
\label{eq:absurd}
\end{align}
But since $\sum_{n = 1}^N \lim_{t^* \rightarrow + \infty} \vartheta_n(t^*)$ $= 1-\alpha > 0$ and $\min_n \{ \bar D_{N+1,n}(t^*) \} > 0$,~\eqref{eq:absurd} can not hold, thus completing the proof. \hfill $\Box$

Requiring the existence of an average node-to-destination multihop path is tantamount to having a Markov transition matrix with a unique absorbing state (the sink node) corresponding to the average graph [cf.~\eqref{eq:deliverability}]. If a node $U_n$ is able to receive packets, but cannot forward them to any other node due to a persistent activity of PU nodes in its proximity (which violates the condition of Proposition 3), then the constraints $t_{i \rightarrow n} = 0$ for all $i \in \cN_{\rightarrow n}$ should be added in (P1). \\

\section{Numerical results}
\label{sec:results}

\begin{figure}
	\centering
  \subfigure[]{\includegraphics[width=0.45\textwidth]{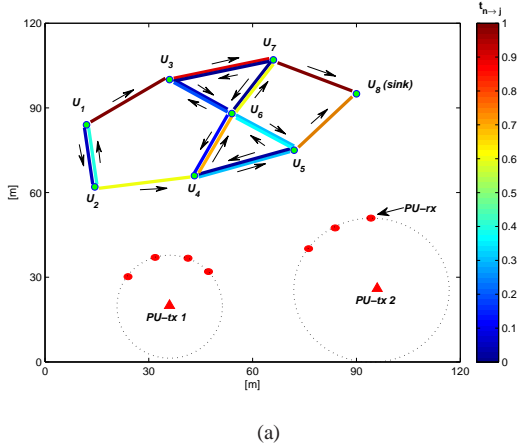}}
  \subfigure[]{\includegraphics[width=0.45\textwidth]{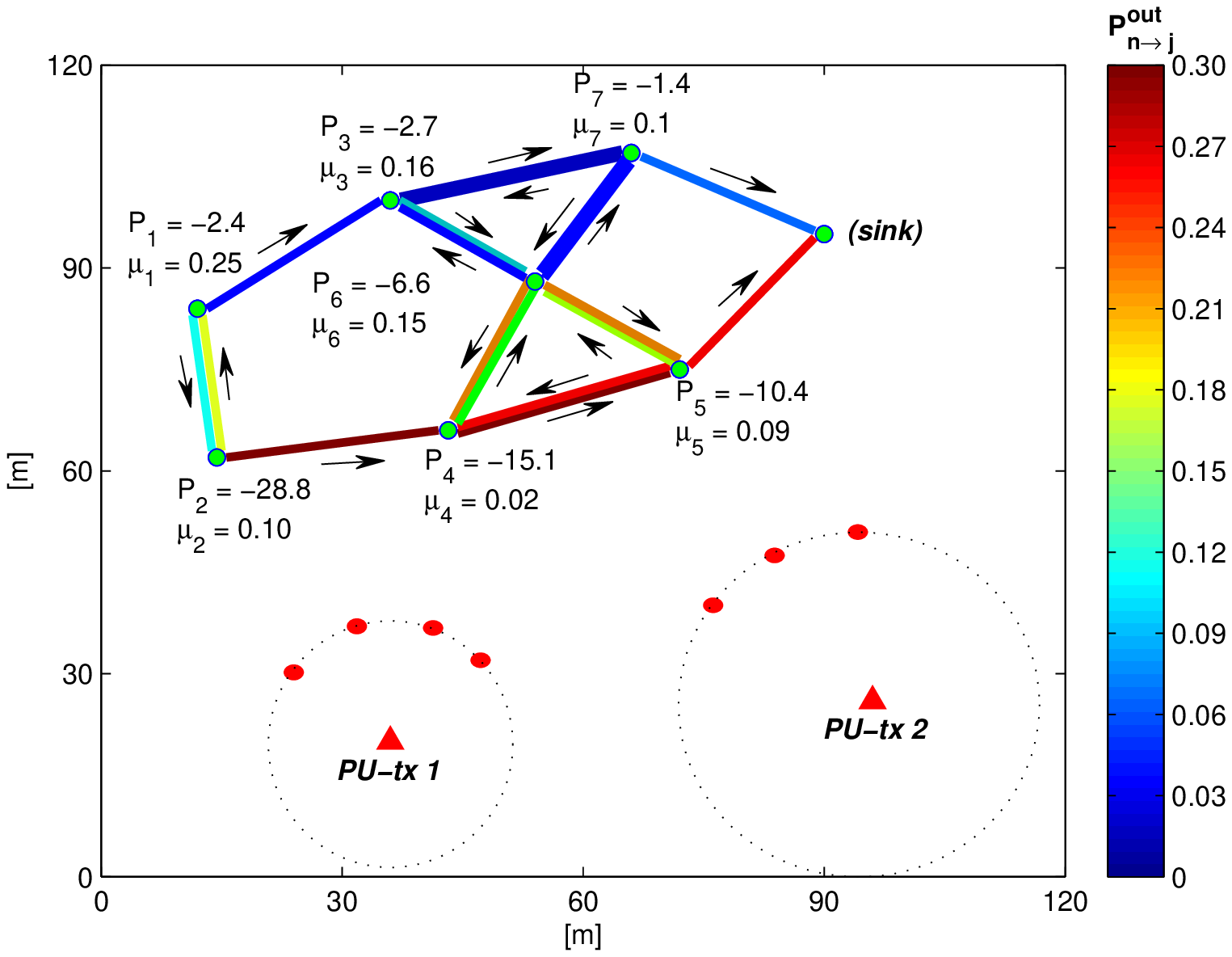}}
  \caption{Test case 1: routing probabilities $\{t_{\ntoi}\}$ (top); and fading-induced outage probabilities (bottom).}
  \label{fig:scenario}
\end{figure}

Consider the scenario depicted in Fig.~\ref{fig:scenario}, where $N=7$ CR nodes cooperate in routing packets to the destination $U_8$.  Two PU sources also transmit to their intended receivers with power
$10$ dBW. In order to protect the PU system without knowing the locations of the PU receivers, $7$ points on the boundary of the PUs' coverage regions are selected~\cite{DKGtvt11}. The PU interference threshold is set to $-80$ dBW. The path loss obeys the model $\|\bx_n - \bx_j\|^{-\eta}$, with $\eta = 3.5$. Log-normal shadowing is generated with standard deviation $6$ dB, and $m = 1$ is used for the small-scale fading (Rayleigh)~\cite{Stu01}. The maximum transmit-power of the CR nodes is set to $\Pnmax = 0$ dBW, and the noise power is $10^{-8}$ W. The SINR threshold $\bar \Gamma_n = -10$ dB, and the sum of exogenous rates $\sum_{n=1}^{N} \rho_n$ is maximized; that is $\cU_n(\rho_n) = \rho_n$ and $\cC_n(P_n)$ = 0, for all $n = 1,\ldots, N$. A larger scale network could also be considered, but the conclusion that one could draw do not depend on the network size. 

Fig.~\ref{fig:scenario}(a) depicts the optimal routing probabilities $\{t_{\ntoi}\}$, obtained by solving (P2) with Algorithm~2. At the first iteration $\ell =1$, a feasible starting point is obtained by properly modifying the approach of~\cite{AvWi71} to the problem at hand, and setting the step-size in~\eqref{eq:projection} equal to $1$. It can be seen that there is a tendency not to route packets through the ``southern'' region of the network; i.e., through nodes that are closer to the PU systems. For example, packets generated by $U_2$ are more likely to be routed through links $U_4 \rightarrow U_6$ and $U_6 \rightarrow U_7$, rather than choosing the shortest path $U_2 \rightarrow U_4 \rightarrow U_5 \rightarrow U_8$. Furthermore, node $U_5$ may decide to send packets to $U_6$ rather than attempting direct transmission to $U_8$ with considerably high probability. This is due to the fact that links starting from and ending to $U_4$ and $U_5$ are characterized by a higher fading- and interference-induced outage probability, as showed in Fig.~\ref{fig:scenario}(b). In fact, not only  PU interference has a detrimental effect on the CR SINRs, but also $U_2$, $U_4$, and $U_5$ are confined to use a lower transmit-power in order to enforce protection of the PU receivers. 
Notice also that $U_2$ may decide to transmit to $U_1$ instead of $U_4$ with considerably high probability. On the other hand, packets generated by $U_1$ and $U_3$ are routed through $U_7$ with high probability, which in this case coincides also with the shortest path. Interestingly, it is necessary to use the primal decomposition algorithm only during the first $5$-$6$ iterations out of the total $14$ (on average) in the successive convex approximation algorithm. In fact, the per-CR interference levels quickly stabilize around steady-state values, with subtle variations for $\ell > 6$.    

\begin{figure}
	\centering
  \subfigure[]{\includegraphics[width=0.45\textwidth]{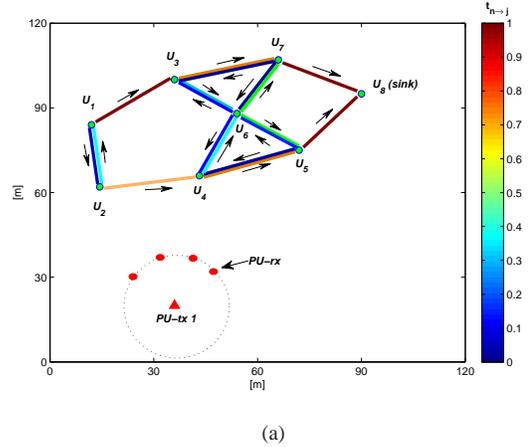}}
  \subfigure[]{\includegraphics[width=0.45\textwidth]{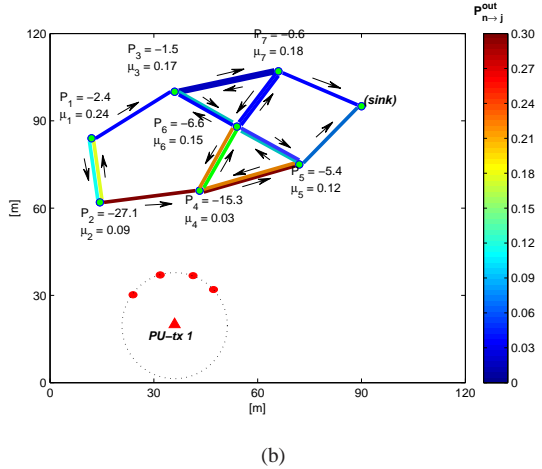}}
  \caption{Test case 2 with the second PU transmitter inactive: routing probabilities $\{t_{\ntoi}\}$ (top); and fading-induced outage probabilities (bottom).}
  \label{fig:scenario2}
\end{figure}

To verify adaptability of the routing probabilities and link reliabilities to the states of the PU systems, consider the case of Fig.~\ref{fig:scenario2}(a), where the same CR network operates only with PU~1 present. Compared to Fig.~\ref{fig:scenario}(a), $U_4$ now  forwards an increased amount of traffic through node $U_5$. As PU~2 is inactive, the outage probability of link $U_4 \rightarrow U_5$ is lower in this case, as confirmed by Fig.~\ref{fig:scenario2}(b). Furthermore, $U_5$ can raise its transmit-power of $10$ dB, which significantly decreases the outage probability of link $U_5 \rightarrow U_8$. As a result, almost none of the packets ($2\%$) are sent to $U_6$. Finally, notice that CR $U_6$ now splits its traffic evenly between $P_5$ and $P_7$. The average exogenous traffic rates, averaged over $20$ different experiments, are reported in Table~\ref{tab:traffic}. It can be seen that $\rho_2$, $\rho_4$, and $\rho_5$ increase in this case. This example demonstrates the capability of the proposed routing approach to adapt routes and transmit-powers to locations of active PUs. 

\begin{figure}
\begin{center}
\includegraphics[width=0.5\textwidth]{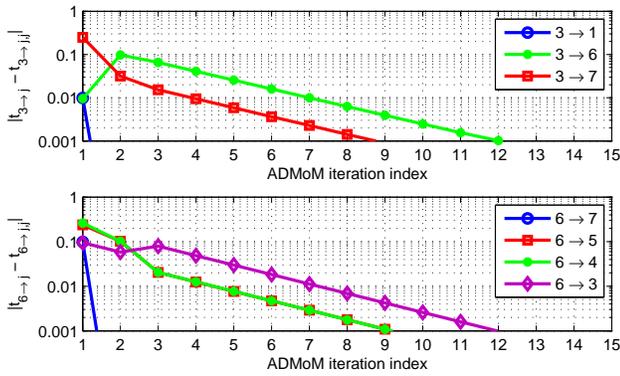}
\caption{Convergence of Algorithm~1.}  
\label{fig:convergence}
\end{center}
\end{figure}

Fig.~\ref{fig:convergence} corroborates the convergence of Algorithm~1 for $\beta = 0.1$ and $c \in \{1, 10\}$. Specifically, the depicted evolution of $|t_{n\rightarrow j}(k) - t_{n\rightarrow j,j}(k)|$ for nodes $U_3$ and $U_6$ shows that the local routing probabilities approximately coincide with those of the neighboring nodes after a few iterations. For example, a gap smaller than $1 \%$ is obtained after $8$ iterations. A similar trend was observed for the transmit-probabilities, which suggests that an online implementation of the algorithm is feasible, and queues will be stable after just a few iterations. 

\begin{table}[h]
\caption{Exogenous traffic rates.}
\begin{center}
\begin{tabular}{l | c | c | c | c | c | c | c }
&  $U_1$ & $U_2$ & $U_3$ & $U_4$ & $U_5$ & $U_6$ & $U_7$ \\
\hline 
Test 1	&  $0.052$ & $0.023$ & $0.028$ & $0.01$ & $0.011$ & $0.008$ & $0.044$ \\
Test 2	&  $0.057$ & $0.038$ & $0.06$ & $0.05$ & $0.015$ & $0.011$ & $0.051$ \\
\hline
\end{tabular}
\end{center}
\label{tab:traffic}
\end{table}%


\section{Conclusions}
\label{sec:conc}
 
A novel cross-layer optimization framework was introduced in this paper. Based on channel and interference level statistics, and the situational awareness provided by spectrum sensing schemes, the novel approach yields optimal routes, transmission probabilities, and transmit-powers. The relevant optimization problem turned out to be non-convex and hence difficult to solve even in a centralized setup. Nevertheless, a successive convex approximation was pursued to find a KKT solution. Primal decomposition and AD-MoM were employed to derive a distributed algorithm, suitable for large networks, and amenable to online implementation. As packets are randomly routed through the network, their deliverability in case of time-varying routing strategies and link reliabilities was asserted. Finally, numerical tests verified the ability of the proposed scheme to adapt network operation to the propagation environment.

\bibliographystyle{IEEEtran}
\bibliography{biblio}

\begin{thebibliography}{10}
\providecommand{\url}[1]{#1}
\csname url@samestyle\endcsname
\providecommand{\newblock}{\relax}
\providecommand{\bibinfo}[2]{#2}
\providecommand{\BIBentrySTDinterwordspacing}{\spaceskip=0pt\relax}
\providecommand{\BIBentryALTinterwordstretchfactor}{4}
\providecommand{\BIBentryALTinterwordspacing}{\spaceskip=\fontdimen2\font plus
\BIBentryALTinterwordstretchfactor\fontdimen3\font minus
  \fontdimen4\font\relax}
\providecommand{\BIBforeignlanguage}[2]{{%
\expandafter\ifx\csname l@#1\endcsname\relax
\typeout{** WARNING: IEEEtran.bst: No hyphenation pattern has been}%
\typeout{** loaded for the language `#1'. Using the pattern for}%
\typeout{** the default language instead.}%
\else
\language=\csname l@#1\endcsname
\fi
#2}}
\providecommand{\BIBdecl}{\relax}
\BIBdecl

\bibitem{Silvester93}
J.~L. Wang and J.~A. Silvester, ``Maximum number of independent paths and radio
  connectivity,'' \emph{IEEE Trans. Wireless Commun.}, vol.~41, no.~10, pp.
  1482--1493, Oct. 1993.

\bibitem{Haenggi05}
M.~Haenggi, ``On routing in random {R}ayleigh fading networks,'' \emph{IEEE
  Trans. Wireless Commun.}, vol.~4, pp. 1553--1562, Jul. 2005.

\bibitem{DeCouto03}
D.~D. Couto, D.~Aguayo, J.~Bicket, and R.~Morris, ``A high-throughput path
  metric for multi-hop wireless routing,'' in \emph{Proc. Int. {ACM} Conf.
  Mobile Computing, Networking}, San Diego, CA, Sep. 2003, pp. 134--156.

\bibitem{ZLZhang07}
H.~Liu, Z.-L. Zhang, J.~Srivastava, and V.~Firoiu, ``Pwave: A multi-source
  multi-sink anycast routing framework for wireless sensor networks,'' in
  \emph{Proc. of Intl. Conf. on Networking}, Atlanta, GA, May 2007.

\bibitem{Ribeiro-infocom07}
A.~Ribeiro, Z.-Q. Luo, N.~Sidiropoulos, and G.~B. Giannakis, ``Modelling and
  optimization of stochastic routing for wireless multihop networks,'' in
  \emph{Proc. IEEE Int. Conf. on Computer Commun.}, Anchorage, AK, May 2007,
  pp. 1748--1756.

\bibitem{Ephremides-twc02}
A.~Ephremides, ``Energy concerns in wireless networks,'' \emph{IEEE Trans.
  Wireless Commun.}, vol.~9, no.~4, pp. 48--59, Aug. 2002.

\bibitem{ZhSa07}
Q.~Zhao and B.~M. Sadler, ``A survey of dynamic spectrum access,'' \emph{IEEE
  Signal Processing Magazine}, vol.~24, no.~3, pp. 79--89, May 2007.

\bibitem{Khalife08}
H.~Khalife, S.~Ahuja, N.~Malouch, and M.~M.~Krunz, ``Probabilistic path
  selection in opportunistic cognitive radio networks,'' in \emph{Proc. of IEEE
  Glob. Telecom. Conf.}, New Orleans, LO, Dec. 2010.

\bibitem{Xin08}
C.~Xin, L.~Ma, and C.-C. Shen, ``A path-centric channel assignment framework
  for cognitive radio wireless networks,'' \emph{Mobile Net. Appl.}, vol.~13,
  no.~5, pp. 463--476, Oct. 2008.

\bibitem{Pefkianakis08}
I.~Pefkianakis, S.~Wong, and S.~Lu, ``{SAMER}: Spectrum aware mesh routing in
  cognitive radio networks,'' in \emph{Proc. of IEEE DySPAN}, Chicago, IL, Oct.
  2008.

\bibitem{abbagnale-Secon}
A.~Abbagnale and F.~Cuomo, ``Connectivity-driven routing for cognitive radio
  ad-hoc networks,'' in \emph{Proceedings of IEEE SECON}, Boston, MA, 2010.

\bibitem{Chowdhury11}
K.~R. Chowdhury and I.~F. Akyildiz, ``{CRP}: A routing protocol for cognitive
  radio ad hoc networks,'' \emph{IEEE J. Sel. Areas Commun.}, vol.~29, no.~4,
  pp. 794--802, Apr. 2011.

\bibitem{KDGjstsp11}
S.-J. Kim, E.~Dall'Anese, and G.~B. Giannakis, ``Cooperative spectrum sensing
  for cognitive radios using {K}riged {K}alman filtering,'' \emph{IEEE J. Sel.
  Topics Sig. Proc.}, vol.~5, pp. 24--36, Feb. 2011.

\bibitem{DBGphycom}
E.~Dall'Anese, J.~A. Bazerque, and G.~B. Giannakis, ``Group sparse {L}asso for
  cognitive network sensing robust to model uncertainties and outliers,''
  \emph{Elsevier Physical Communication}, Nov. 2011.

\bibitem{Lott06}
C.~Lott and D.~Teneketzis, ``Stochastic routing in ad-hoc networks,''
  \emph{IEEE Trans. Auto. Contr.}, vol.~51, no.~1, pp. 52--70, Jan. 2006.

\bibitem{Stu01}
G.~L. St\"uber, \emph{Principles of Mobile Communication}, 2nd~ed.\hskip 1em
  plus 0.5em minus 0.4em\relax Boston, MA: Kluwer Academic Publishers, 2001.

\bibitem{DKGGPtwc11}
E.~Dall'Anese, S.-J. Kim, G.~B. Giannakis, and S.~Pupolin, ``Power control for
  cognitive radio networks under channel uncertainty,'' \emph{IEEE Trans.
  Wireless Commun.}, vol.~10, pp. 3541--3551, Dec. 2011.

\bibitem{Fen60}
L.~F. Fenton, ``The sum of lognormal probability distributions in scatter
  transmission systems,'' \emph{IRE Trans. Commun. Syst.}, vol.~8, no.~1, pp.
  57--67, Mar. 1960.

\bibitem{Shin04}
S.~Shin, S.~Choi, H.~S. Park, and W.~H. Kwon, ``Packet error rate analysis of
  {IEEE} 802.15.4 under {IEEE} 802.11b interference,'' in \emph{Proc. of WWIC
  2005, LNCS, Springer}, May. 2004, pp. 279--288.

\bibitem{Beaulieu94}
A.~A. Abu-Dayya and N.~C. Beaulieu, ``Comparison of methods of computing
  correlated lognormal sum distributions and outages for digital wireless
  applications,'' in \emph{Proc. IEEE Veh. Tech. Conf.}, May 1994, pp.
  175--179.

\bibitem{Raotit88}
R.~Rao and A.~Ephremides, ``On the stability of interacting queues in a
  multi-access system,'' \emph{IEEE Trans. Info. Theory}, vol.~34, pp.
  918--930, Sep. 1988.

\bibitem{Georgiadis06}
M.~J.~N. L.~Georgiadis and L.~Tassiulas, ``Resource allocation and cross-layer
  control in wireless networks,'' \emph{Found. Trends in Netw.}, vol.~1, no.~1,
  pp. 1--144, 2006.

\bibitem{Loynes}
R.~Loynes, ``The stability of a queue with non-independent interarrival and
  service times,'' \emph{Mathematical Proc. of the Cambridge Philosophical
  Society}, vol.~58, pp. 497--520, 1962.

\bibitem{DKGtvt11}
E.~Dall'Anese, S.-J. Kim, and G.~B. Giannakis, ``Channel gain map tracking via
  distributed {K}riging,'' \emph{IEEE Trans. Veh. Technol.}, vol.~60, no.~3,
  pp. 1205--1211, Mar. 2011.

\bibitem{Zhang09}
R.~Zhang, ``On peak versus average interference power constraints for
  protecting primary users in cognitive radio networks,'' \emph{IEEE Trans.
  Wireless Commun.}, vol.~8, no.~4, pp. 2112--2120, Apr. 2009.

\bibitem{Marques_book}
A.~G. Marques, N.~Gatsis, and G.~B. Giannakis, ``Optimal cross-layer design of
  wireless fading multi-hop networks,'' \emph{in Cross Layer Designs in WLAN
  Systems}, {N.} Zorba, C. Skianis, and C. Verikoukis, Eds. Leicester, UK:
  Troubador Pub., 2011.

\bibitem{Ribeiro}
A.~Ribeiro, ``Wireless cooperative communications and networking,'' {PhD}
  thesis, University of Minnesota. 2006. [Online]:
  \texttt{http://www.seas.upenn.edu/}$\sim$\texttt{aribeiro/preprints/}.

\bibitem{Chianitwc03}
M.~Chiani, D.~Dardari, and M.~K. Simon, ``New exponential bounds and
  approximations for the computation of error probability in fading channels,''
  \emph{IEEE Trans. Wireless Commun.}, vol.~2, no.~4, pp. 840--845, Jul. 2003.

\bibitem{Haggman04}
N.~Ermolova and S.-G. Haggman, ``Simplified bounds for the complementary error
  function; application to the performance evaluation of signal processing
  systems,'' in \emph{Proc. of the 12th European Signal Proces. Conf.}, Vienna,
  Austria, Sep. 2004.

\bibitem{MaWr78}
B.~R. Marks and G.~P. Wright, ``A general inner approximation algorithm for
  nonconvex mathematical programs,'' \emph{Oper. Res.}, vol.~26, no.~4, pp.
  681--683, Jul.-Aug. 1978.

\bibitem{BoV04}
S.~Boyd and L.~Vandenberghe, \emph{Convex Optimization}.\hskip 1em plus 0.5em
  minus 0.4em\relax Cambridge University Press, 2004.

\bibitem{Ribeiro-twc08}
A.~Ribeiro, N.~Sidiropoulos, and G.~B. Giannakis, ``Optimal distributed
  stochastic routing algorithms for wireless multihop networks,'' \emph{IEEE
  Trans. Wireless Commun.}, vol.~7, no.~11, pp. 4261--4272, Nov. 2008.

\bibitem{Palomar06}
D.~P. Palomar and M.~Chiang, ``A tutorial on decomposition methods for network
  utility maximization,'' \emph{IEEE J. Sel. Areas Commun.}, vol.~24, no.~9,
  pp. 1439--1451, 2006.

\bibitem{Nedic09}
A.~Nedi\'{c} and A.~Ozdaglar, ``Approximate primal solutions and rate analysis
  for dual subgradient methods,'' \emph{SIAM J. Optim.}, vol.~19, no.~4, pp.
  1757--1780, 2009.

\bibitem{BeT89}
D.~P. Bertsekas and J.~N. Tsitsiklis, \emph{Parallel and Distributed
  Computation: Numerical Methods}.\hskip 1em plus 0.5em minus 0.4em\relax
  Englewood Cliffs, NJ: Prentice-Hall, 1989.

\bibitem{Michelot86}
C.~Michelot, ``A finite algorithm for finding the projection of a point onto
  the canonical simplex of $\mathbb{R}^n$,'' \emph{J. Optim. Theory Appl.},
  vol.~50, no.~1, pp. 195--200, Jul. 1986.

\bibitem{Newman}
M.~E.~J. Newman, \emph{Networks: An Introduction}.\hskip 1em plus 0.5em minus
  0.4em\relax Oxford University Press, 2010.

\bibitem{AvWi71}
M.~Avriel and A.~C. Williams, ``An extension of geometric programming with
  applications in engineering optimization,'' \emph{Journal of Engineering
  Mathematics}, vol.~5, no.~3, pp. 187--194, Jul. 1971.

\end{thebibliography}

\begin{IEEEbiography}[{\includegraphics[width=1in,height=1.25in,clip,keepaspectratio]{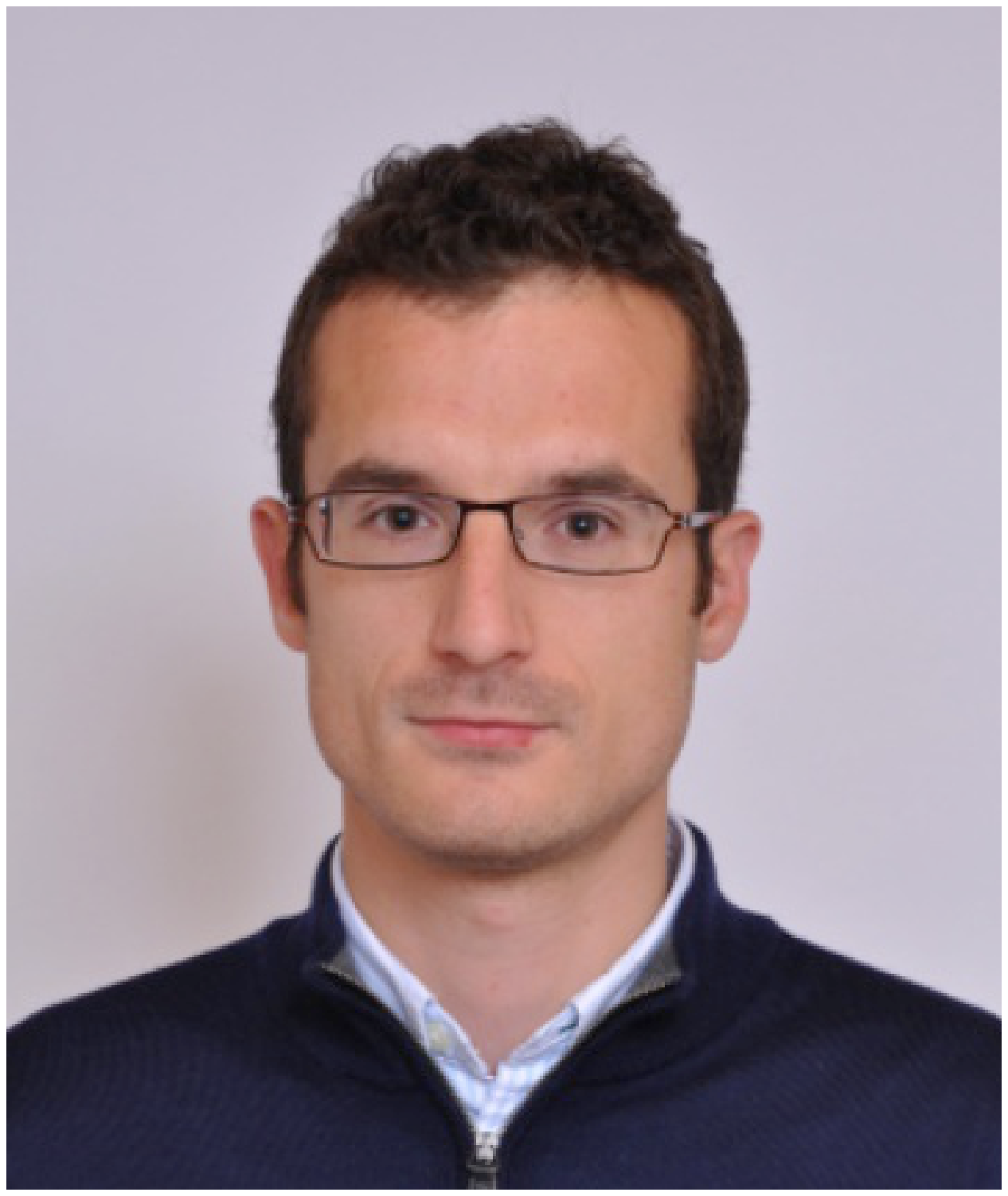}}]{Emiliano Dall'Anese (S'08, M'11)}
received the Laurea Triennale (B.Sc degree) and the Laurea Specialistica (M.Sc degree) in Telecommunications Engineering from the University of Padova, Italy, in 2005 and 2007, respectively, and the Ph.D in Information Engineering at the Department of Information Engineering (DEI), University of Padova, Italy, in 2011. From January 2009 to September 2010 he was a visiting scholar at the Department of Electrical and Computer Engineering, University of Minnesota, USA. He is currently a post-doctoral associate at the Department of Electrical and Computer Engineering, University of Minnesota, USA. 

His research interests lie in the areas of statistical signal processing, communication theory, and networking. Current research focuses on wireless cognitive radio systems, IP networks, and power distribution networks.
\end{IEEEbiography}

\newpage

\begin{IEEEbiography}[{\includegraphics[width=1in,height=1.25in,clip,keepaspectratio]{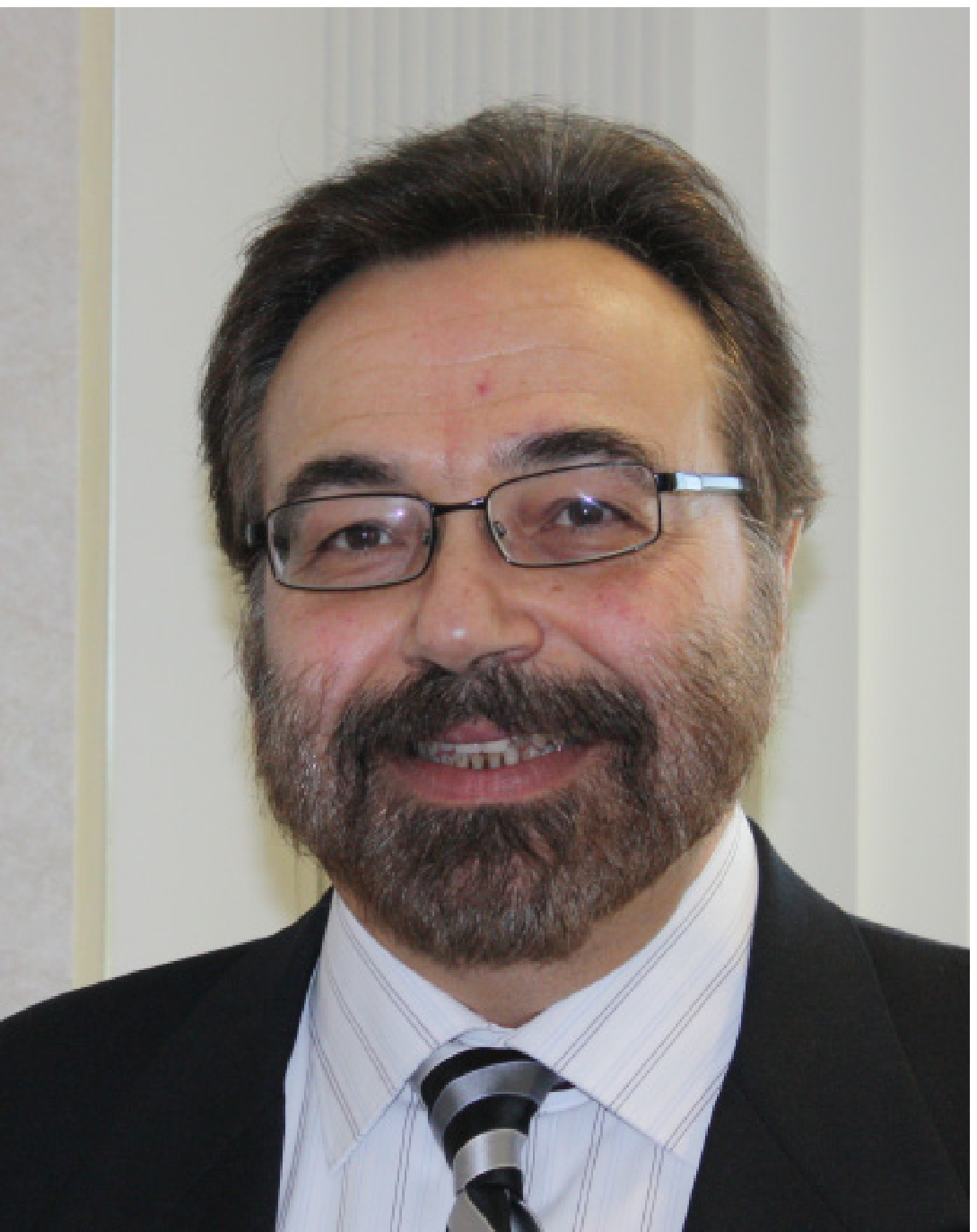}}]{Georgios B. Giannakis (F'97)}
received his Diploma in Electrical 
Engr. from the Ntl. Tech. Univ. of Athens, Greece, 1981. From 
1982 to 1986 he was with the Univ. of Southern California (USC), 
where he received his MSc. in Electrical Engineering, 1983, MSc. 
in Mathematics, 1986, and Ph.D. in Electrical Engr., 1986. Since 
1999 he has been a professor with the Univ. of Minnesota, where 
he now holds an ADC Chair in Wireless Telecommunications in the 
ECE Department, and serves as director of the Digital Technology
Center. 

His general interests span the areas of communications, networking 
and statistical signal processing - subjects on which he has 
published more than 300 journal papers, 500 conference papers, 
20 book chapters, two edited books and two research monographs. 
Current research focuses on compressive sensing, cognitive radios, 
network coding, cross-layer designs, wireless sensors, social and 
power grid networks. He is the (co-) inventor of twenty patents 
issued, and the (co-) recipient of eight paper awards from the IEEE 
Signal Processing (SP) and Communications Societies, including the 
G. Marconi Prize Paper Award in Wireless Communications. He also 
received Technical Achievement Awards from the SP Society (2000), from 
EURASIP (2005), a Young Faculty Teaching Award, and the G. W. Taylor 
Award for Distinguished Research from the University of Minnesota. 
He is a Fellow of EURASIP, and has served the IEEE in a number of posts, 
including that of a Distinguished Lecturer for the IEEE-SP Society. 
\end{IEEEbiography}

\end{document}